\documentclass[12pt]{article}
\begin{document}
\begin{center}
{\bf On periodic boundary value problem for the Sturm-Liouville operator}
\\
\medskip\medskip
{Alexander Makin}
\end{center}

We consider the eigenvalue problem for the Sturm-Liouville operator
$$
Lu=u''-q(x)u \eqno(1)
$$
with boundary conditions
$$
u(0)\mp u(1)=0, \quad u'(0)\mp u'(1)=0 \eqno(2)
$$
where either the sign "$-$" is chosen two times (case 1) or the sign "+" is chosen two times (case 2),
i.e. boundary conditions (2) are periodic or antiperiodic. Function $q(x)$ is an
arbitrary complex-valued function of the class $L_1(0,1)$.

Let $\{u_n(x)\}$ be the system of eigenfunctions and associated functions of problem (1)+(2).
It is well known that this system is complete and minimal in $L_2(0,1)$. We denote
$$
\alpha_n=\int_0^1q(x)e^{2\pi inx}dx, \quad \beta_n=\int_0^1q(x)e^{-2\pi inx}dx.
$$
Suppose the function $q(x)$ satisfies the following conditions: 
$q(x)\in W_1^m[0,1]$, $q^{(j)}(0)=q^{(j)}(1)$ where $j=\overline{0,m-1}$, $m=0,1,\ldots$.

{\bf Theorem 1.} {\it If for all even (in the case 1) or odd (in the case 2)
$n>n_0$ where $n_0$ is a natural number
$$
|\alpha_n|>\frac{c_0}{n^{m+1}}, \quad 0<c_1<|\frac{\alpha_n}{\beta_n}|<c_2
$$
$(c_0>0)$, then the root function system $\{u_n(x)\}$ of corresponding problem (1)+(2)
forms a Riesz basis for $L_2(0,1)$.}

{\bf Theorem 2.} {\it If there exists a sequence of even (in the case 1) or 
odd (in the case 2) numbers $n_k$ $(k=1,2,\ldots)$ such that 
$$
|\alpha_{n_k}|>\frac{c_0}{n_k^{m+1}}, \quad 
|\beta_{n_k}|>\frac{c_0}{n_k^{m+1}} 
$$
$(c_0>0)$ moreover 
$
\lim_{k\to\infty}(|\alpha_{n_k}/\beta_{n_k}|+|\beta_{n_k}/\alpha_{n_k}|)=\infty
$
, then the root function system $\{u_n(x)\}$
of corresponding problem (1)+(2)
is not a basis for $L_2(0,1)$.}

It is easy to verify that the function
$$
q(x)=\sum_{n=1}^{\infty}\gamma_n\left(\frac{e^{2\pi inx}}{n^{\varepsilon_1}}+\frac{e^{-2\pi inx}}{n^{\varepsilon_2}}\right)
$$
satisfies all conditions of Theorem 2. Here $0<\varepsilon_1<\varepsilon_2<1$ and also in the case 1 $\gamma_n=1 $ if
$n=2^p$ and $\gamma_n=0$ if $\gamma_n\ne2^p$, and in the case 2 $\gamma_n=1$ if
$n=2^p+1$ and $\gamma_n=0$ if $n\ne2^p+1$ $(p=1,2, \ldots)$.

We denote by $Q$ the set of potentials $q(x)$ such that the system
$\{u_n(x)\}$ is a Riesz basis for $L_2(0,1)$, $\bar Q=L_1(0,1)\setminus Q$. 
From Theorem 1 and Theorem 2 it is easy to obtain the following

{\bf Corollary.} {\it The sets $Q$ and $\bar Q$ are dense everywhere in $L_1(0,1)$.}

Convergence of spectral expansions corresponding to problem (1)+(2) was studied
by O.A. Veliev, N.B. Kerimov, Kh.P. Mamedov.

\medskip
Moscow State Academy of Instrument-Making and Informatics, Stromynka 20, Moscow,
107996, Russia

\medskip
E-mail address: alexmakin@yandex.ru

\end{document}